\documentclass[11pt]{article}
\usepackage{graphicx}
\usepackage{amssymb,amsfonts,amsthm}
\usepackage{amsmath,amsthm,amssymb}
\usepackage{amsfonts}

\newtheorem{theorem}{Theorem}[section]
\newtheorem{lemma}[theorem]{Lemma}

\theoremstyle{definition}

\newtheorem{rk}[theorem]{Remark}

\newcommand{\la}{\langle}
\newcommand{\ra}{\rangle}
\newcommand{\Con}{\mathrm{Con}}
\newcommand{\dist}{{\mathrm{dist}}}

\newcommand{\Lab}{\phi}

\newcommand{\ttt}{{\cal T}}

\newcommand{\aaa}{{\cal A}}

\newcommand{\vk}{van Kampen }

\newcommand{\iv}{^{-1}}

\begin{document}

\renewcommand{\theequation}{\thesection.\arabic{equation}}

\title{Groups with non-simply connected asymptotic cones}
\author{A.Yu. Ol'shanskii, M.V. Sapir\thanks{Both authors were supported in part by
the NSF grant DMS 0245600. In addition, the research of the first
author was supported in part by the Russian Fund for Basic
Research 02-01-00170 and by the Russian Fund for the Support of
Leading Research Groups 1910.2003.1}}
\date{}
\maketitle

\begin{abstract} We construct a group (an HNN extension of a free group)
with polynomial isoperimetric function, linear isodiametric function
and non-simply connected asymptotic cones.
\end{abstract}

\section{Introduction}

Recall that a function $f\colon \mathbb{N}\to \mathbb N$ is an
isoperimetric (resp. isodiametric) function  of a finite
presentation $\la X\mid R\ra$ of a group $G$ if every word $w$ in
$X$, that is equal to $1$ in $G$, is freely equal to a product of
conjugates $\prod_{i=1}^m x_i\iv r_i x_i$ where $r_i$ or $r_i\iv$ is
in $R$, $x_i\in (X\cup X\iv)^*$, and $m\le f(|w|)$ (resp. $|x_i|\le
f(|w|)$ for every $i$).

Isoperimetric (resp. isodiametric) functions $f_1$, $f_2$ of any two
finite presentations of the same group $G$ are {\em equivalent},
that is $f_2(n)\le Cf_1(Cn) +Cn$, $f_1(n)\le Cf_2(Cn)+Cn$ for some
constant $C$. As usual, we do not distinguish equivalent functions.
The smallest isoperimetric function of a group is called its {\em
Dehn function}.

In terms of \vk diagrams, $f$ is an isoperimetric (isodiametric)
function of the finite presentation if for every word $w$, that is
equal to $1$ in $G$ there exists a \vk diagram with boundary label
$w$ and area (resp. diameter) at most $f(|w|)$ (see Gersten
\cite{Gersten} for details).

Let $(X,\dist)$ be a metric space, $o=(o_n)$ be a sequence of points
in $X$, $d=(d_n)$ be an increasing sequence of numbers with $\lim
d_n=\infty$, and let $\omega\colon P(\mathbb{N})\to \{0,1\}$ be an
ultrafilter. An {\em asymptotic cone} $\Con_\omega(X,o,d)$ of
$(X,\dist)$ is the subset of the cartesian power $X^{\mathbb N}$
consisting of sequences $(x_i)_{i\in \mathbb N}$ with
$\limsup\frac{\dist(o_i,x_i)}{d_i}<\infty$ where we identify two
sequences $(x_i)$ and $(y_i)$ with $\lim_\omega
\frac{\dist(x_i,y_i)}{d_i}=0$. The distance between two elements
$(x_i)$ and $(y_i)$ in the asymptotic cone $\Con_\omega(X,o,d)$ is
defined as $\lim_\omega \frac{\dist(x_i,y_i)}{d_i}$. Here
$\lim_\omega$ is the $\omega$-limit defined as follows. If $a_n$ is
a bounded sequence of real numbers then $\lim_\omega(a_n)$ is the
(unique) number $a$ such that for every $\epsilon>0$,
$\omega(\{n\mid |a_n-a|<\epsilon\})=1$.

If $G$ is a finitely generated group, then {\em asymptotic cones} of
$G$ are, by definition, asymptotic cones of the Cayley graph of $G$
(with respect to some generating set). Asymptotic cones of $G$ do
not depend on the choice of the sequence $o$, so we can always
assume that $o=(1)$ where $1$ is the identity, and use notation
$\Con_\omega(G,d)$.

Gromov proved \cite{Gromov} that if all asymptotic cones of a group
$G$ (for all $\omega$ and all $d$) are simply connected then $G$ is
finitely presented, has polynomial isoperimetric function and linear
isodiametric function. Papasoglu \cite{Papasoglu} proved that if a
finitely presented group has quadratic isoperimetric function then
all its asymptotic cones are simply connected. But in general the
existence of polynomial isoperimetric functions does not imply that
the asymptotic cones are simply connected. Indeed, it is not even
true that polynomial isoperimetric inequality implies linear
isodiametric inequality. First examples of groups with polynomial
Dehn functions and non-linear isodiametric functions were
constructed in \cite{Bridson} and \cite{SBR}\footnote{Paper
\cite{SBR} appeared as a preprint in 1997, two years earlier than
\cite{Bridson}. But results of \cite{Bridson} were announced four
years earlier, in 1993.}. The question of whether we can guarantee
that the asymptotic cones are simply connected by requiring that
both the Dehn function is polynomial and the isodiametric function
is linear, was open. The question was mentioned, in particular, by
Dru\c tu in \cite{Drutu}. The following theorem answers this
question.

Let $$G=\la\theta_1,\theta_2, a, k\mid a^{\theta_i}=a,
k^{\theta_i}=ka, i=1,2\ra$$ where $x^y=y\iv xy$. Thus $G$ has a
balanced presentation with $4$ generators and $4$ relators. It is
clear that $G$ is a split extension of the free group $\la a,k\ra$
by the free group $\la\theta_1, \theta_2\ra$.

\begin{theorem}\label{thmain} The group $G$ has a cubic
isoperimetric function, a linear isodiametric function, and no
simply connected asymptotic cones.
\end{theorem}

This group is an $S$-machine in terminology of \cite{OSlogn} or a
hub-free interpretation of an $S$-machine in terminology of
\cite{SBR}. As an $S$-machine, $G$ has one state letter $k$, one
tape letter $a$, and two rules, both of the form $[k\to ka]$.
Applying results from \cite{SBR} to $G$, we can deduce that $G$ has
a cubic isoperimetric function (in fact its Dehn function is exactly
$n^3$ by \cite{OSlogn}), and linear isodiametric function. For the
sake of completeness, we present below a direct proof of these
statements. Then we shall prove that $G$ has no simply connected
asymptotic cones.

\begin{rk} Using the construction from \cite{OSlogn}, we can
find an $S$-machine (a multiple HNN extension of a free group) with
Dehn function $n^2\log n$, linear isodiametric function and no
simply connected asymptotic cones. (One can also assume that the
conjugacy problem for that group is undecidable.) The proof is only
technically more difficult than the proof of Theorem \ref{thmain},
and the group is much more complicated, so we do not include this
example here leaving it as an advanced exercise.

It would be interesting to find out what the asymptotic cones of $G$
are. In particular, how many non-by-Lipschitz-equivalent asymptotic
cones does $G$ have, and what are the fundamental groups of these
cones? More generally, it is interesting to find out what are the
asymptotic cones of $S$-machines? Topological properties of
asymptotic cones  may reflect computational properties of the
$S$-machines (for a general definition of an $S$-machine see
\cite{SBR} or \cite{OSlogn}).
\end{rk}

\section{Proof}

For every letter $x$, an $x$-edge in a \vk diagram is an edge
labeled by $x^{\pm 1}$; an $x$-cell is a cell whose boundary
contains an $x$-edge. If $x$ is a free letter in an HNN extension,
then an $x$-band in a diagram over $G$ is a sequence of cells
containing $x$-edges, such that every two consecutive cells share an
$x$-edge. Our group $G$ can be considered as an HNN extension with
free letters $\theta_1, \theta_2$ (obviously), and also as an
HNN-extension with free letter $k$. Thus we can consider
$\theta$-bands (i.e. $\theta_i$-bands, $i=1,2$) and $k$-bands in \vk
diagrams. It is also convenient to consider $a$-bands, that is
sequences of cells corresponding to the relations $a^{\theta_i}=a$,
such that every two consecutive cells share an $a$-edge. The
boundary of the union of cells from an $x$-band $\ttt$ has the form
$e\iv pfq\iv$ where $e, f$ are the only $x$-edges on the boundary.
The paths $p$ and $q$ are called the {\em sides} of the band $\ttt$,
$e$ and $f$ are called the {\em start} and {\em end} edges of the
band. The {\em median} of the band $\ttt$ is a polygonal simple line
that connects the midpoints of $e, f$ and is contained in the
interior of the cells composing the band. We fix one median for each
band. A median of an $x$-annulus is defined similarly; it is a
simple closed curve.

If $x, y$ are two different letters in $\{\theta,k,a\}$ then one can
consider an $(x,y)$-annulus that is a union of an $x$-band and a
$y$-band sharing the first and the last cells only, and such that
the region bounded by the medians of the bands does not contain the
start and end edges of the band.

We shall say that a \vk diagram over $G$ is reduced if it does not
contain a pair of cells that share a boundary edge and are mirror
images of each other.

The following lemma is a particular case of Lemma 3.3 from
\cite{OSlogn}.

\begin{lemma}\label{NoAnnul}
A reduced van Kampen diagram $\Delta$ over $G$ has no $k$-annuli,
$\theta$-annuli, $(k,\theta)$-annuli,  $a$-annuli,
$(a,\theta)$-annuli.
\end{lemma}
\proof We assume that $\Delta$ is a counterexample with minimal
area. This means in particular that  the boundary of $\Delta$ is the
boundary component of an annulus $\aaa$, where $\aaa$ has one of the
types from the formulation of the lemma.

(1) Let $\aaa$ be a $k$-annulus. Then it consists of
$(k,\theta)$-cells. Hence there is a maximal $\theta$-band $\ttt$ in
$\Delta$, whose first cell $\pi_1$ and the last cell $\pi_2$ belong
to $\aaa$. Being members of the same $k$-band $\aaa$ and
$\theta$-band $\ttt$, the cells $\pi_1$ and $\pi_2$ cannot be
neighbors in $\ttt$ (the diagram is reduced). Hence $\ttt$ and a
part of $\aaa$ form a $(\theta,q)$-annulus.  The area of the
subdiagram bounded by this annulus is smaller than that of $\Delta$.
This contradicts the choice of $\Delta$.

(2) Let $\aaa$ be a $\theta$-annulus. If it contains $k$-cells, then
we come to a contradiction as in (1). Otherwise $\Delta$ has no
$k$-cells since there is no counter-example of smaller area. The
inner part of $\aaa$ (i.e. the subdiagram bounded by the median of
$\aaa$) has no $\theta$-edges for the same reason. Hence $\Gamma$
has no cells corresponding to the relations of the group $G$. So, on
the one hand, the inner label of $\aaa$ is a cyclically reduced
non-empty word in $Y$ since $\Delta$ is a reduced diagram, and on
the other hand, this word is freely equal to 1, a contradiction.

(3) Let $\aaa$ be a $(k,\theta)$-annulus. Then the maximal $k$-band
$\ttt$ of $\aaa$ cannot have more than two cells because otherwise
$\Delta$ would contain a smaller counterexample as in (2). Hence the
length of $\ttt$ is 2, and its cells are mirror copies of each other
(since they belong to the same $\theta$-band), a contradiction (we
assumed that the diagram is reduced).

(4) Let $\aaa$ be an $a$-annulus. Then its boundary labels are words
in $\theta_1, \theta_2$. This leads, as in (1), to a smaller
$(a,\theta)$-annulus, a contradiction.

(5) Let $\aaa$ be a $(\theta,a)$-annulus and let $\ttt$ be the
maximal $a$-band of it. It cannot have more than two
$(\theta,a)$-cells because otherwise there would be a smaller
$(\theta,a)$-annulus. Hence the length of $\ttt$ is 2, and its cells
are mirror copies of each other, a contradiction. The lemma is
proved.
\endproof

The following lemma is a part of Theorem \ref{thmain}.

\begin{lemma} \label{cubic} $G$ has a cubic isoperimetric function and a linear
isodiametric function.
\end{lemma}

\proof Let $\Delta$ be a reduced \vk diagram over $G$ with perimeter
$n$. We need to estimate the area of $\Delta$  and its diameter in
terms of $n$.

Suppose that $\Delta$ is reduced. Every $k$-cell in $\Delta$ is an
intersection of a $\theta$-band and a $k$-band. Since the bands
intersect only once by Lemma \ref{NoAnnul}, the number of $k$-cells
does not exceed the product of the number of maximal $\theta$-bands
and the number of maximal $k$-bands. By Lemma \ref{NoAnnul}, every
maximal $\theta$-band and every maximal $k$-band connect two edges
on the boundary of $\Delta$. Hence the number of maximal $k$- and
$\theta$-bands does not exceed $n$, and the number of $k$-cells does
not exceed $\frac{n^2}{4}$. Every $a$-cell that is not a $k$-cell is
the intersection of a $\theta$-band and an $a$-band. Every $a$-band
in $\Delta$ either starts on the boundary of $\Delta$ or on the
boundary of a $k$-cell. Since each $k$-cell contains exactly one
$a$-edge, the number of maximal $a$-bands is at most
$\frac{n^2}{8}+\frac{n}2$. Since an $a$-band and a $\theta$-band can
intersect only once (Lemma \ref{NoAnnul} again) the total number of
$a$-cells that are not $k$-cells in $\Delta$ is a most
$(\frac{n^2}{8}+\frac{n}{2})\frac{n}2=\frac{n^3}{16}+\frac{n^2}{4}$,
and the total number of cells in $\Delta$ does not exceed
$\frac{n^3}{16}+\frac{n^2}{2}$.

Let $s$ be the number of $\theta$-bands in $\Delta$. Then $s\le
\frac{n}2$. By Lemma \ref{NoAnnul}, there exists a $\theta$-band
that whose side is a part of the boundary $\partial\Delta$. Then
$\Delta$ is obtained by gluing $\ttt$ and a reduced diagram
$\Delta_1$ with $s-1$ $\theta$-bands. Every vertex on a side of
$\ttt$ can be connected to the boundary of $\Delta$ by a path of
length at most 2. Therefore by induction on $s$, we can deduce
that every vertex inside $\Delta$ can be connected to the boundary
of $\Delta$ by a path of length at most $2s\le 2\frac{n}2=n$.
Hence the diameter of $\Delta$ is at most $5n/2$.
\endproof

\begin{rk} 1. As we have mentioned before, the Dehn function of $G$ is
$n^3$. To obtain the lower bound, it is enough to consider the
diagrams with boundary label $[k^n, \theta_1^n\theta_2^{-n}]$ where
$[x,y]=x\iv y\iv xy$ (see below). One can also use the general
statement from \cite{OSlogn}.

2. The second part of the proof of Lemma \ref{cubic} works without
any significant change for all multiple HNN extensions of free
groups. Thus the isodiametric function of any such group is linear.
\end{rk}

It remains to show that every asymptotic cone of $G$ is not simply
connected. Suppose that an asymptotic cone $\Con_\omega(G,d)$ is
simply connected. It was noticed by Gromov \cite{Gromov} (see also
\cite{Papasoglu} or \cite{Drutu} for more details) that then for
every $M>1$ there exists a number $k$ such that for every constant
$C\ge 1$, every loop $l$  in the Cayley complex of $G$, such that
$\frac{1}{C}d_m\le |l|\le Cd_m$ for any sufficiently large $m$,
bounds a (singular) disc that can be subdivided into $k$ (singular)
subdiscs with perimeters at most $\frac{|l|}{M}$. (In fact one can
assume that the statement is true for {\em some} $M>1$; then one can
deduce the statement for all $M$ by further subdividing the
subdiscs.)

Fix a natural number $n=d_m$, $m>>1$. Let $u_n$ be the commutator
$[k^n, \theta_1^n\theta_2^{-n}].$ Clearly, $u_n=1$ in $G$ since
$k^{\theta_1^n}=k^{\theta_2^n}=ka^n$. The corresponding \vk
diagram has the form of a trapezium \cite{SBR}, \cite{OSlogn} with
the top and the bottom sides $p, p'$ labeled by $k^n$ and the left
and right sides $q, q'$, labeled by $\theta_1^n\theta_2^{-n}$. The
perimeter of that diagram is $6n$. Let us show that the loop in
the Cayley graph of $G$ corresponding to $u_n$ cannot bound a disc
decomposed into at most $l\le \sqrt{n}$ subdiscs of perimeter $n$.
That would contradict the statement from the previous paragraph
for $M=C=6$. Suppose that such a decomposition exists. Then we
have a (not necessarily reduced) \vk diagram $\Delta$ with
boundary label $u_n$ composed of $l$ subdiagrams
$\Delta_1,...,\Delta_l$ such that the perimeter of each $\Delta_i$
is at most $n$. Consider any $\theta_1$-edge $e$ of the path $q$.
The $\theta_1$-band in $\Delta$ that starts at $e$ cannot end on
$p, p'$ (since these paths do not contain $\theta$-edges), or on
$q$ (since every edge of $q$ points from the initial vertex $q_-$
to the terminal vertex $q_+$ of $q$). Hence it ends on $q'$. Since
$\theta$-bands do not intersect, they connect corresponding
$\theta$-edges on $q$ and $q'$, that is the $\theta$-edge number
$i$ on $q$  is connected with the $\theta$-edge number $i$ on
$q'$.

Let $e$ be the $\theta_1$-edge number $n$ on $q$ and $\ttt$ be the
maximal $\theta$-band starting at $e$. Let $r$ be the top side of
$\ttt$. Then the label $\Lab(r)$ belongs to the free group $\la k,
a\ra$ and is equal to $\theta_1^{-n}k^n\theta_1^n$ in $G$. Hence
this word is freely equal to $(ka^n)^n$. Since the number of
subdiagrams $\Delta_i$ is $l\le \sqrt{n}$, there is a subpath $w$
of $r$ such that the initial and terminal vertices of $w$ belong
to the boundary of one of the $\Delta_i$ and the freely reduced
form $W$ of $\Lab(w)$ contains $ka^nk$ as a subword. Hence, in the
group $G$, we should have an equality $W=U$ where $|U|\le
\frac{n}2$ since the perimeter of $\Delta_i$ does not exceed $n$.

Let $\Gamma$ be a reduced diagram for this equality with boundary
$r_1r_2\iv$ where $\Lab(r_1)=W$, $\Lab(r_2)=U$. Let $r_1'$ be a
subpath of $r_1$ with label $ka^nk$. Consider the two maximal
$k$-bands $\ttt_1$ and $\ttt_2$ starting on the $k$-edges $e$ and
$f$ of $r_1'$.

Note that since $r_1$ does not contain $\theta$-edges, every
$\theta$-band crossing $\ttt_1$ or $\ttt_2$ must end on $r_2$.
Hence if $\ttt_1$ or $\ttt_2$ end on $r_1$ or the length of that
$k$-band must be zero. But this would mean that a non-trivial
subword of a freely reduced word $W$ is equal to 1 in the free
group $\la k, a\ra$ which is impossible. Hence $\ttt_1$ and
$\ttt_2$ end on $r_2$

Let $\Gamma'$ be the subdiagram of $\Gamma$ bounded by
$r_1''=r_1'\setminus\{e,f\}$, sides of $\ttt_1$ and $\ttt_2$, and
a part $\bar r_2$ of $r_2$. Let $\partial\Gamma' = \bar r_1\bar
r_2\iv=r_1''(r_2')\iv$. The length of a side of $\ttt_1$ or
$\ttt_2$ is at most twice as large as the number of $\theta$-edges
in it, and the number of $\theta$-edges in these sides is at most
$|r_2|-|\bar r_2|-2$ since there are two $k$-edges on $r_2$. Hence
we have $|r_2'| < 2|r_2|-2|\bar r_2|+|\bar r_2| = 2|r_2|-|\bar
r_2|$.

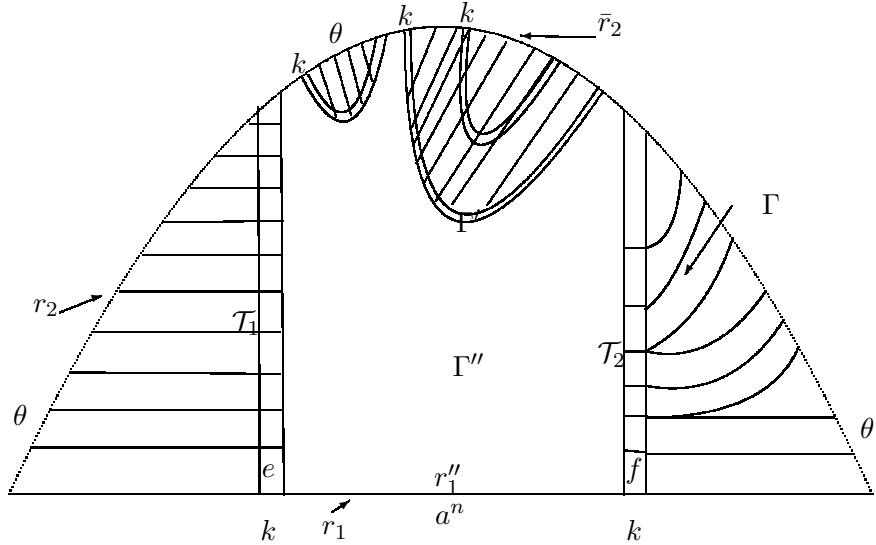
\unitlength .8mm 
\linethickness{0.4pt}
\ifx\plotpoint\undefined\newsavebox{\plotpoint}\fi 
\begin{picture}(145.75,102.75)(0,0)
\put(2.25,12.25){\line(1,0){143.5}}
\qbezier(2.25,12.25)(74.5,167.75)(145.75,12.25)
\multiput(44,12.25)(-.04167,10.75){6}{\line(0,1){10.75}}
\multiput(48,12)(-.04167,5.60417){12}{\line(0,1){5.60417}}
\put(104.5,12){\line(0,1){64.25}}
\put(108.25,12.25){\line(0,1){60}}
\put(45.5,6.25){\makebox(0,0)[cc]{$k$}}
\put(45.5,16.25){\makebox(0,0)[cc]{$e$}}
\put(106.25,6.25){\makebox(0,0)[cc]{$k$}}
\put(106.25,16.25){\makebox(0,0)[cc]{$f$}}
\put(75.75,9){\makebox(0,0)[cc]{$a^n$}}
\qbezier(51,81.5)(60,64)(65,88.5)
\qbezier(52,82)(59.38,67)(63.25,88)
\qbezier(68,89)(68.5,31.38)(101,79.25)
\qbezier(69,89.25)(69.38,33.38)(100.25,80)
\qbezier(77.75,89.75)(74.25,54.13)(93.75,84)
\qbezier(78.75,89.25)(75.25,57.5)(92.75,84.75)
\put(78.75,33.75){\makebox(0,0)[cc]{$\Gamma''$}}
\put(50.5,84){\makebox(0,0)[cc]{$k$}}
\put(68.25,92){\makebox(0,0)[cc]{$k$}}
\put(78.5,92.5){\makebox(0,0)[cc]{$k$}}
\put(75.25,14.75){\makebox(0,0)[cc]{$r_1''$}}
\put(129,60.5){\makebox(0,0)[cc]{$\Gamma$}}
\put(114.75,48.75){\vector(-2,-3){.09}}\multiput(122.5,60.25)(-.04211957,-.0625){184}{\line(0,-1){.0625}}
\put(78.75,57.75){\makebox(0,0)[cc]{$\Gamma'$}}
\put(56.5,6.25){\makebox(0,0)[cc]{$r_1$}}
\put(59,11.25){\vector(3,2){.09}}\multiput(56,9.25)(.0625,.0416667){48}{\line(1,0){.0625}}
\put(8.5,43.25){\makebox(0,0)[cc]{$r_2$}}
\put(17.75,45.25){\vector(3,1){.09}}\multiput(10.75,42.5)(.1060606,.0416667){66}{\line(1,0){.1060606}}
\put(41.75,40.75){\makebox(0,0)[cc]{$\ttt_1$}}
\put(102.5,35.5){\makebox(0,0)[cc]{$\ttt_2$}}
\put(102.25,90.25){\makebox(0,0)[cc]{$\bar r_2$}}
\put(87.5,88.5){\vector(-1,0){.09}}\multiput(99.75,88.25)(-2.04167,.04167){6}{\line(-1,0){2.04167}}
\put(6,20){\line(1,0){41.75}}
\put(9,26.25){\line(1,0){38.75}}
\multiput(12.5,32.5)(5.83333,-.04167){6}{\line(1,0){5.83333}}
\put(16.25,39.25){\line(1,0){31.25}}
\put(20.75,46){\line(1,0){26.75}}
\put(24.5,52){\line(1,0){23.25}}
\multiput(28,57.5)(3.33333,.04167){6}{\line(1,0){3.33333}}
\put(32.5,63.25){\line(1,0){14.75}}
\put(36.75,68.5){\line(1,0){10.75}}
\put(42.25,73.75){\line(1,0){5}}
\put(108.5,19){\line(1,0){34}}
\put(108.5,25){\line(1,0){31}}
\qbezier(108.5,25)(128.5,25.13)(133.5,36.75)
\qbezier(108.25,30.25)(122.38,27.75)(131,41.25)
\qbezier(108.25,36)(119.38,33.5)(128,46)
\qbezier(108.5,36.25)(118,41.5)(122.5,54.75)
\qbezier(108,43)(113,46.5)(118,61)
\qbezier(108.25,53.25)(113.38,54.38)(114,66)
\multiput(54,83.5)(.0416667,-.1287879){66}{\line(0,-1){.1287879}}
\multiput(56.5,85)(.0416667,-.1477273){66}{\line(0,-1){.1477273}}
\multiput(58.75,86)(.0416667,-.162037){54}{\line(0,-1){.162037}}
\multiput(61,87)(.0416667,-.1428571){42}{\line(0,-1){.1428571}}
\multiput(73.25,89.75)(-.0420354,-.1039823){113}{\line(0,-1){.1039823}}
\put(68.5,78){\line(0,1){0}}
\multiput(76.75,89.75)(-.04213483,-.09269663){178}{\line(0,-1){.09269663}}
\multiput(81.5,88.5)(-.042124542,-.077838828){273}{\line(0,-1){.077838828}}
\multiput(78.75,89.5)(-.04204545,-.0875){220}{\line(0,-1){.0875}}
\multiput(85,87.5)(-.04204893,-.071865443){327}{\line(0,-1){.071865443}}
\multiput(89.75,86.5)(-.042091837,-.06505102){392}{\line(0,-1){.06505102}}
\multiput(92.75,85)(-.042085427,-.061557789){398}{\line(0,-1){.061557789}}
\multiput(96.75,82)(-.042134831,-.061095506){356}{\line(0,-1){.061095506}}
\put(145,23.25){\makebox(0,0)[cc]{$\theta$}}
\put(4.25,25.5){\makebox(0,0)[cc]{$\theta$}}
\put(56.75,89){\makebox(0,0)[cc]{$\theta$}}
\put(104.75,53.25){\line(1,0){3.75}}
\put(104.5,43.5){\line(1,0){3.75}}
\put(104.5,36){\line(1,0){4}}
\put(104.5,30.25){\line(1,0){3.75}}
\put(104.5,25.25){\line(1,0){4}}
\multiput(104.5,19.5)(.625,-.04167){6}{\line(1,0){.625}}
\end{picture}

Every
$k$-band $\ttt$ in $\Gamma'$ connects two edges $e(\ttt), f(\ttt)$ on $\bar r_2$. We
say that a $k$-band $\ttt$ is {\em farther} from the boundary than a
$k$-band $\ttt'$ if the subpath of $\bar r_2$ between $e(\ttt)$,
$f(\ttt)$ contains $e(\ttt')$ and $f(\ttt')$. Note that the length
of a side of $\ttt$ is at most twice as large as the length of the
subpath of $\bar r_2$ between $e(\ttt)$ and $f(\ttt)$.

Let us remove from $\Gamma'$ all the $k$-bands that are the
farthest from $\bar r_2$  together with parts of $\Gamma'$ between
their sides and $\bar r_2$. Let us denote the resulting subdiagram
of $\Gamma'$ by $\Gamma''$. Then $\partial\Gamma''=r_1''r_2''$.
Since $k$-bands do not intersect, $|r_2''|\le |r_2'|+|\bar r_2|$.
Hence $|r_2''|< 2|r_2|-|\bar r_2|+|\bar r_2|=2|r_2|\le n$.

Note that $\Gamma''$ does not contain $k$-cells. Therefore every
$a$-band that starts on $r_1''$ must end on $r_2''$. But that is
impossible since $|r_2''|<n$, and $\Lab(r_1'')=a^n$, a
contradiction. This completes the proof of Theorem \ref{thmain}.

\begin{minipage}[t]{3 in}
\noindent Alexander Yu. Ol'shanskii\\ Department of Mathematics\\
Vanderbilt University \\ alexander.olshanskiy@vanderbilt.edu\\
http://www.math.vanderbilt.edu/$\sim$olsh\\ and\\ Department of
Higher Algebra, MEHMAT\\
Moscow State University\\
olshan@shabol.math.msu.su\\
\end{minipage}
\begin{minipage}[t]{3 in}
\noindent Mark V. Sapir\\ Department of Mathematics\\
Vanderbilt University\\
m.sapir@vanderbilt.edu\\
http://www.math.vanderbilt.edu/$\sim$msapir\\
\end{minipage}

\addtocontents{toc}{\contentsline {section}{\numberline {
}References \hbox {}}{\pageref{bibbb}}}

\begin{thebibliography}{1}
\label{bibbb}

\bibitem{Bridson} Martin Bridson. Asymptotic cones and polynomial
isoperimetric inequalities. Topology 38 (1999), no. 3, 543--554.

\bibitem{Drutu} Cornelia Dru\c tu. Quasi-isometry invariants and
asymptotic cones. International Conference on Geometric and
Combinatorial Methods in Group Theory and Semigroup Theory (Lincoln,
NE, 2000). Internat. J. Algebra Comput. 12 (2002), no. 1-2, 99--135.

\bibitem{Gersten} Steve M. Gersten. Isoperimetric and isodiametric functions
of finite presentations. Geometric group theory, Vol. 1 (Sussex,
1991), 79--96, London Math. Soc. Lecture Note Ser., 181, Cambridge
Univ. Press, Cambridge, 1993.

\bibitem{Gromov} M. Gromov.
Asymptotic invariants of infinite groups. Geometric group theory,
Vol. 2 (Sussex, 1991), 1--295, London Math. Soc. Lecture Note Ser.,
182, Cambridge Univ. Press, Cambridge, 1993.

\bibitem{OSlogn} A.Yu. Olshanskii and M. V. Sapir. Groups with small Dehn functions and
bipartite chord diagrams. preprit, arXiv math.GR/0411174.

\bibitem{Papasoglu} P. Papasoglou. Asymptotic cones and the quadratic isoperimetric
inequality.  Journal of Differential Geometry, 1996, 44, 789--806.

\bibitem{SBR} M. V. Sapir, J. C. Birget, E. Rips.
\newblock Isoperimetric and isodiametric functions of groups,
Annals of Mathematics, 157, 2 (2002), 345-466.
\end{thebibliography}
\end{document}